\documentclass[12pt]{article}
\usepackage{mathrsfs}
\usepackage{amsfonts}
\usepackage{latexsym, amssymb, amsmath, amscd, amsfonts, epsfig, graphicx, colordvi}
\usepackage{latexsym,amsmath,amssymb,amsfonts,epsfig,graphicx,cite,psfrag}
\usepackage{eepic,color,colordvi,amscd,amsthm}
\usepackage{pstcol,pstricks,color}

\newtheorem{theorem}{Theorem}[section]

\newtheorem{lem}[theorem]{Lemma}

\def\qed{\hfill \rule{4pt}{7pt}}

\def\pf{\noindent {\it Proof.} }

\newcommand{\inv}{\rm{inv}}

\begin{document}

\begin{center}
{\Large\bf Labeled Partitions with Colored Permutations}
\end{center}

\begin{center}
William Y.C. Chen$^1$, Henry Y. Gao$^2$, Jia He$^3$

Center for Combinatorics, LPMC-TJKLC\\
Nankai University, Tianjin 300071, P.R. China

$^1$chen@nankai.edu.cn, $^2$gaoyong@cfc.nankai.edu.cn,
$^3$hejia1@msu.edu

\end{center}

\vskip 6mm \noindent{\small {\bf Abstract.} In this paper, we extend
the notion of labeled partitions with ordinary permutations to
colored permutations in the sense that the colors are endowed with a
cyclic structure. We use labeled partitions with colored
permutations to derive the generating function  of the
$\mathrm{fmaj}_k$ indices of colored  permutations. The second
result is a combinatorial treatment of a relation on the
$q$-derangement numbers with respect to colored permutations which
leads to the formula of Chow for signed permutations and the formula
of Faliharimalala and Zeng \cite{Fal082} on colored permutations.
The third result is an involution on permutations that implies the
generating function formula for the signed $q$-counting of the major
indices due to Gessel and Simon. This involution can be extended to
signed permutations. In this way, we obtain a combinatorial
interpretation of a formula of Adin, Gessel and Roichman.



\section{Introduction }

In this paper, we will be concerned with the wreath product $S_n^k=
C_k\wr S_n$ of the symmetric group on $[n]=\{1, 2, \ldots, n\}$ and
the cyclic group $C_k$ on $\{0, 1, \ldots, k-1\}$ whose elements are
considered being arranged on a cycle, see Adin and Roichman
\cite{Adin01} and Wagner \cite{Wag03}. The elements in $S_n^k$ are
also called colored permutations \cite{Bag07}. The derangements with
respect to the group $S_n^k$ are studied by Faliharimalala and Zeng
\cite{Fal08, Fal082}.

We will extend the notion of labeled partitions with ordinary
permutations to colored permutations.  A $k$-colored permutation is
written as $\pi(1)_{c_1}\, \pi(2)_{c_2}\, \cdots\, \pi(n)_{c_n}$,
where $\pi(1)\, \pi(2)\, \cdots\,  \pi(n)$ is a permutation on $[n]$
and $c_i\in \{0, 1, \ldots, k-1\}$.  For example, $4_2\, 3_0\, 1_2
\, 5_0 \, 2_1$ is a colored permutation in $S_5^3$. We  define a
total order on the elements of $S_n^k$ as follows
\begin{equation}
1_{k-1}<2_{k-1}<\cdots<n_{k-1}<1_{k-2}<2_{k-2}<\cdots<n_{k-2}<\cdots
<1_0<2_0<\cdots<n_0.\label{Order}
\end{equation}
We now recall the following definitions:
\begin{align}
D(\sigma)&:=\{i\in[n-1]:\sigma(i)>\sigma(i+1)\},\nonumber\\[5pt]
\mathrm{maj}(\sigma)&:=\sum_{i\in D(\sigma)}i,\nonumber\\[5pt]
N_j(\sigma)&:=\#\{i\in[n]:\sigma(i)\ \text{has subscript $j$}
\},\quad j=1,\ldots,k-1,\nonumber\\[5pt]
\mathrm{fmaj_k}(\sigma)&:=k\mathrm{maj}(\sigma)+N_1(\sigma)+2N_2(\sigma)+\cdots
+(k-1)N_{k-1}(\sigma).\label{flagmajor}
\end{align}

The set $D(\sigma)$ is called the descent set of $\sigma\in S_n^k$.
It should be noted that  Adin and Roichman  \cite{Adin01} give the
definition of flag major index of an element in $S_n^k$ by the
unique factorization into Coxeter elements, and they prove that
$\mathrm{fmaj}_k$ has the above expression (\ref{flagmajor}). In
this paper, we will consider the formula (\ref{flagmajor}) as the
definition of the $\mathrm{fmaj}_k$ index. From this point of view,
our approach may be regarded as purely combinatorial.

For $k=1$, $S_n^1$ is usually written as $S_n$.  For $k=2$, $S_n^2$
becomes the group of signed permutations on $[n]$, often denoted by
$B_n$, and the minus sign is often denoted by a bar. Moreover, the
$\mathrm{fmaj}_k$ index reduces to the $\mathrm{fmaj}$ index for
signed permutations as defined by
\[ \mathrm{fmaj}(\pi)=2\mathrm{maj}(\pi)+N(\pi),\]
where $N(\pi)$ denotes the number of  negative  elements of $\pi$
and $\mathrm{maj}(\pi)$ is defined with respect to the following
order \[ \bar{1}< \bar{2}<\cdots < \bar{n} < 1 < 2< \cdots < n.\]

Using labeled partitions with colored permutations, we get the
generating function of the $\mathrm{fmaj_k}$ indices on $S_n^k$,
\begin{align}
 \sum_{\pi\in S_n^k}q^{\mathrm{fmaj}_k(\pi)}=[k]_q[2k]_q\cdots
[nk]_q,\label{eq generating function of S_n^3}
\end{align}
where $[k]_q=1+q+q^2+\cdots +q^{k-1}$. The above formula is a
natural extension of the formulas for the generating functions for
the major index and the fmaj index, see Faliharimalala and Zeng
\cite{Fal082}. Recall that for the cases of ordinary permutations
and signed permutations we have
\begin{align}
\sum_{\pi \in S_n}q^{\mathrm{maj}(\pi)}=[n]!\label{gsn}
\end{align}
and
\begin{align}
\sum_{\pi \in
B_n}q^{\mathrm{fmaj}(\pi)}=[2]_q[4]_q\cdots[2n]_q\label{gbn}.
\end{align}

 The second result is
a combinatorial treatment of a relation on the $q$-derangement
numbers $D_n^k(q)$ with respect to $S_n^k$.  This relation implies
the formula for $d_n^k(q)$ by the $q$-binomial inversion, as given
by Faliharimalala and Zeng \cite{Fal082}.  For $n\geq 1,$ let
\[ \mathscr{D}_n:=\{\sigma\in S_n:\sigma(i)\neq i\ \text{for all}\ i\in
[n]\}\]  be the set of  derangements on $S_n$. Gessel defined the
$q$-derangement  numbers  by
\[ d_n(q):=\sum_{\sigma\in\mathscr{D}_n}q^{\mathrm{maj}(\sigma)}\]
and proved  that
\begin{align}
d_n(q)=[n]_q!\sum_{k=0}^n\frac{(-1)^kq^{k \choose
2}}{[k]_q!},\label{eq1}
\end{align}
where $[n]_q!=[1]_q[2]_q\cdots [n]_q.$ Wachs \cite{Wachs89} found a
combinatorial proof of the above formula.  Later, Chow \cite{Chow06}
generalized Wachs's argument  to type $B$ derangements. Similarly,
Chow defined
\[ \mathscr{D}_n^B:=\{\sigma\in B_n:\sigma(i)\neq i\ \text{for all}\
i\in [n]\}\] as the set of derangements in $B_n$ and
\[
d_n^B(q):=\sum_{\sigma\in\mathscr{D}_n^B}q^{\mathrm{fmaj}(\sigma)}.\]
 Chow has shown that
\begin{align}
d_n^B(q)=[2]_q[4]_q\cdots[2n]_q\sum_{k=0}^n\frac{(-1)^kq^{2{k\choose
2}}}{[2]_q[4]_q\cdots[2k]_q} \label{eq2}.
\end{align}

 The notion of derangements of type $B$ can be generalized to $S_n^k$, as given by
Faliharimalala and Zeng \cite{Fal082}.  We define \[
\mathscr{D}_n^k:=\{\sigma\in S_n^k:\sigma(i)\neq i_0\ \text{for
all}\ i\in [n]\}\] and
 \[
 d_n^k(q):=\sum_{\sigma\in\mathscr{D}_n^k}q^{\mathrm{fmaj}_k(\sigma)}.\]
Faliharimalala and Zeng have shown that
\begin{align}
d_n^k(q)=[k]_q[2k]_q\cdots[nk]_q\sum_{j=0}^n\frac{(-1)^jq^{k{j\choose
2}}}{[k]_q[2k]_q\cdots[jk]_q}.
\end{align}

The argument of Chow for $d_n^B(q)$ can be extended to the case of
$d_n^k(q)$. Our proof is based on the structure of labeled
partitions with colored permutations, which is an extension of the
combinatorial approach of Chen and Xu \cite{Chen06} for ordinary
permutations. We will present the proof for the case $k=3$, which is
essentially a proof for the general case.

The third result is concerned  the following formula of Gessel and
Simon \cite{Wachs92} on the signed $q$-counting of permutations
 with respect to the major index:
 \[ \sum_{\pi\in
S_n}\mathrm{sign}(\pi)q^{\mathrm{maj}(\pi)}
=[1]_{q}[2]_{-q}[3]_{q}[4]_{-q}\cdots[n]_{(-1)^{n-1}q}.\]

Note that a combinatorial proof of the above formula has been given
by Wachs \cite{Wachs92} based on permutations. We will present an
involution on labeled partitions that leads to a combinatorial
interpretation of the above formula. Moreover, our involution can be
extended to signed permutations. This gives a combinatorial proof of
the following formula of Adin-Gessel-Roichman \cite{Adin05} for the
signed $q$-counting of signed permutations with respect to the
$\mathrm{fmaj}$ index:
\[ \sum_{\pi\in
B_n}\mathrm{sign}(\pi)
q^{\mathrm{fmaj}(\pi)}=[2]_{-q}[4]_{q}\cdots[2n]_{(-1)^{n}q}.\]

\section{Labeled Partitions and the $\mathrm{fmaj}_k$ Index}\label{Labeled Partitions on $S_n^3$}

In this section, we introduce the notion of labeled partitions with
colored permutations.  Using labeled partitions, we give a
combinatorial proof of the following formula for the generating
function of the $\mathrm{fmaj}_k$ indices of colored permutations in
$S_n^k$, given by Haglund,  Loehr and  Remmel \cite{Hag05}, see
also, Faliharimalala and Zeng \cite{Fal082}.

\begin{theorem}\label{theo generating function fmaj_k}
We have
\begin{align*}
 \sum_{\pi\in S^k_n}q^{\mathrm{fmaj}_k(\pi)}=[k]_q[2k]_q\cdots
[nk]_q.
\end{align*}
\end{theorem}

Recall that given a colored permutation $\pi\in S_n^k$, $N_j(\pi)$
denotes the number of elements $\pi(i) \in \pi$ with subscript $j$,
where $j=1, 2, \ldots, k-1$. The $\mathrm{fmaj}_k$ index which is
originally defined algebraically by Adin and Roichman has the
following equivalent form
\begin{equation*} \mathrm{fmaj}_k(\pi) = k \mathrm{maj}(\pi) +
N_1(\pi) + 2N_2(\pi) + \cdots + (k-1)N_{k-1}(\pi).
\end{equation*}

Clearly,  Theorem \ref{theo generating function fmaj_k} is a
generalization of the formulas (\ref{gsn}) and (\ref{gbn}) for
permutations and signed permutations. We now proceed to give a
combinatorial proof of Theorem \ref{theo generating function fmaj_k}
by using labeled partitions with colored permutations.

 Let
$\lambda=(\lambda_1,\lambda_2,\ldots,\lambda_n)$ be an integer
partition with at most $n$ parts where
$\lambda_1\geq\lambda_2\geq\cdots\geq \lambda_n\geq 0$. We adopt the
notation in Andrews  \cite{Andrews76}.  We write
$|\lambda|=\lambda_1+\lambda_2+\cdots+\lambda_n.$ A \textit{labeled
partition} associated with $S_n^3$ is defined as a pair
$(\lambda,\pi)$, where $\lambda$ is a partition with at most $n$
parts and $\pi=\pi(1)\pi(2)\cdots\pi(n)$ is a colored permutation in
$S_n^3$. We can also employ the  two-row notation to represent a
labeled partition
\[\left(
\begin{array}{cccc}
\lambda_1&\lambda_2&\cdots&\lambda_n\\
\pi(1)&\pi(2)&\cdots&\pi(n)
\end{array}\right).
\]

A labeled partition $(\lambda, \pi)$ is said to be \textit{standard}
if $\pi(i)>\pi(i+1)$ implies $\lambda_i>\lambda_{i+1}.$ It is easy
to see that a labeled partition $(\lambda,\pi)$ is standard if
$\lambda_i=\lambda_{i+1}$ implies $\pi(i)<\pi(i+1)$.

Given a colored element $w_i$, we use $c(w_i)$ to denote the color
or subscript $i$, and use $d(w_i)$ to denote the element $w$ after
removing the color $i$.

Let $P_n^3$ denote the set of partitions with at most $n$ parts such
that each part is divisible by $3$. Given $\pi\in S_n^3$, we denote
by $Q_{\pi}$  the set of standard labeled  partitions such that
$\lambda_i-\mathrm{c}(\pi(i))$ is divisible by $3$.

\begin{lem}\label{lem g_pi}
Given $\pi\in S_n^3$, there is a bijection $g_{\pi}\colon
\lambda\rightarrow (\mu, \pi)$ from $ P_n^3$ to
  $Q_{\pi}$ such that $|\lambda|+\mathrm{fmaj}_3(\pi)=|\mu|$.
\end{lem}

\pf We define $\mu$ as follows:
\begin{align*}
\mu=(\lambda_1+3a_1+\mathrm{c}(\pi(1)),\lambda_2+3a_2+\mathrm{c}(\pi(2)),\ldots,\lambda_n+3a_n+\mathrm{c}(\pi(n))),
\end{align*}
where $a_i$ is the number of descents in
$\pi(i)\pi(i+1)\cdots\pi(n)$. From the above definition, it is clear
that $\mu$ is a partition and $\mu_i-\mathrm{c}(\pi(i))$ is
divisible by $3$. We only need to show that $(\mu,\pi)$ is standard.
We have the following cases.

\noindent Case 1: $\lambda_i>\lambda_{i+1}$.  In this case, we have
$\lambda_i+3a_i+\mathrm{c}(\pi(i))=\mu_i>\mu_{i+1}=\lambda_{i+1}+3a_{i+1}+\mathrm{c}(\pi(i+1)),$
since $\lambda_i-\lambda_{i+1}\geq 3, a_i\geq a_{i+1}$ and
$|\mathrm{c}(\pi(i))-\mathrm{c}(\pi(i+1))|<3$.

\noindent  Case 2: $\lambda_i=\lambda_{i+1}$. We further consider
the following two subcases:
\begin{itemize}
\item[(i)] If $\pi(i)>\pi(i+1)$, then it is easy to verify that
\[ \lambda_i+3a_i+\mathrm{c}(\pi(i))=\mu_{i}>\mu_{i+1}
=\lambda_{i+1}+3a_{i+1}+\mathrm{c}(\pi(i+1)).\]
\item[(ii)] If $\pi(i)<\pi(i+1)$ and $\pi(i),\pi(i+1)$ have the
same subscript, then we have
\[ \lambda_i+3a_i+\mathrm{c}(\pi(i))=\mu_{i}
=\mu_{i+1}=\lambda_{i+1}+3a_{i+1}+\mathrm{c}(\pi(i+1)).\]
 Otherwise,
if $\pi(i)$ and $\pi(i+1)$ have different subscripts, then we see
that the subscript of $\pi(i)$ is greater than that of $\pi(i+1)$.
This implies that
\[ \lambda_i+3a_i+\mathrm{c}(\pi(i))=\mu_{i}>\mu_{i+1}
=\lambda_{i+1}+3a_{i+1}+\mathrm{c}(\pi(i+1)).\]
\end{itemize}

Now we see that the labeled partition $(\mu,\pi)$ is standard.
Conversely, given a labeled  partition $(\mu, \pi)\in Q_\pi$ , we
can uniquely recover the  partition $\lambda\in P_n^3$.  \qed

Consequently, we obtain the following  formula.

\begin{theorem}\label{theo fmaj_3 generating function}
For $n\geq 1$, we have
\begin{align*}
 \sum_{\pi\in S_n^3}q^{\mathrm{fmaj}_3(\pi)}=[3]_q[6]_q\cdots
[3n]_q.\label{eq generating function of S_n^3}
\end{align*}
\end{theorem}
\pf We consider the following equivalent form of (\ref{eq generating
function of S_n^3}):
\begin{align*}
\frac{1}{(q^3;q^3)_n}\sum_{\pi \in
S_n^3}q^{\mathrm{fmaj}_3(\pi)}=\frac{1}{(1-q)^n},
\end{align*}
where \[ (q^3;q^3)_n=(1-q^3)(1-q^6)\cdots (1-q^{3n}).\]
 Let $W_n$ be the set of  sequences of $n$ nonnegative integers.
It is clear that  $\frac{1}{(q^3;q^3)_n}$ and $\frac{1}{(1-q)^n}$
are the generating functions for numbers of partitions in $P_n^3$
and $W_n$, respectively. Therefore, it suffices to construct a
bijection $\phi\colon(\lambda\,,\pi)\rightarrow s$ from $(P_n^3,
S_n^3)$ to $W_n$ such that $|\lambda|+\mathrm{fmaj_3}(\pi)=|s|$,
where $|s|$ denotes the sum of entries of $s$. The bijection $\phi$
can be described as follows:

\noindent Step 1. Use the bijection in Lemma \ref{lem g_pi} to
derive a standard labeled partition $(\mu,\pi)$ from
$(\lambda,\pi)$.

\noindent
 Step 2. Based on the two row representation of the labeled partition $(\mu,\pi)$,
           we permute the columns to make the second row become the identity permutation
           by ignoring the subscripts of the elements in $\pi$. Let $s$ denote the first row
           of the array.

It is not difficult to see that the above procedure is reversible.
The inverse of $\phi$ consists of four steps.

\noindent
 Step 1. For a sequence
 $s=(s(1),s(2),\ldots, s(n))\in W_n$,
we construct a two row array
\[\left(
\begin{array}{cccc}
s(1)&s(2)&\cdots&s(n)\\
1&2&\cdots&n
\end{array}\right).
\]

\noindent
 Step 2. For each element $i\in[n]$, we may construct a colored permutation
 $1_{c_1}2_{c_2}\cdots n_{c_n}$, where $c_i= s(i) \pmod{3}$.
 Clearly, we have $s^*(i)=s(i)-c_i$ is divisible by $3$. So we are
 led to the following array
\[\left(
\begin{array}{cccc}
s^*(1)&s^*(2)&\cdots&s^*(n)\\
 1_{c_1}&2_{c_2}&\cdots&n_{c_n}
\end{array}\right).
\]
\noindent Step 3. Permute the columns of the above array to make the
first row  $s^*(j_1)s^*(j_2)\cdots s^*(j_n)$  in decreasing order.
Moreover, we order the elements in the second row in increasing
order if they correspond to the same elements in the first row.  We
denote the resulted labeled partition by
\[\left(
\begin{array}{cccc}
s^*(j_1)&s^*(j_2)&\cdots&s^*(j_n)\\[5pt]
{\delta(1)}_{e_{1}} &{\delta(2)}_{e_{2}}&\cdots&{\delta(n)}_{e_{n}}
\end{array}\right).
\]

\noindent
 Step 4.  Recover the initial labeled partition
$(\lambda,\pi)$
 from  the array produced in  Step
$3$ by the following rule:
\[(\lambda^{*},\pi)=\left(
\begin{array}{cccc}
s^*(j_1)-3a_1&s^*(j_2)-3a_2&\cdots&s^*(j_n)-3a_n\\[5pt]
{\delta(1)}_{e_{1}} &{\delta(2)}_{e_{2}}&\cdots&{\delta(n)}_{e_{n}}
\end{array}\right),
\]
 where $a_k$ is the number of descents in
${\delta(k)}_{e_{k}}\cdots{\delta(n)}_{e_{n}}$.

It is easy to see that the above procedure is feasible. Moreover,
one can verify that $\phi\cdot \phi^{-1}=\mathrm{id}$ and
$\phi^{-1}\cdot \phi=\mathrm{id}$, where $\mathrm{id}$ denotes the
identity map. This completes the proof. \qed

Let us give an example. Let $n=7$, $\lambda=(18,18,18,9,9,6,3)$ and
$\pi=3_2\,4_2\,6_0\,5_1\,7_2\,2_1\,1_2$. Then we obtain
$s=(5,10,29,29,16,27,14)$ via the following steps:
\begin{align*}\left(
\begin{array}{ccccccc}
18&18&18&9&9&6&3\\
3_2&4_2&6_0&5_1&7_2&2_1&1_2
\end{array}\right)&\stackrel{\mathrm{Step \;1}}{\longrightarrow}
\left(
\begin{array}{ccccccc}
29&29&27&16&14&10&5\\
3_2&4_2&6_0&5_1&7_2&2_1&1_2
\end{array}\right)\\[6pt]
&\stackrel{\mathrm{Step
\;2}}{\longrightarrow}\,(5,10,29,29,16,27,14).\end{align*} The
reverse process from $s$ to $(\lambda, \pi)$ are demonstrated as
follows:
\begin{align*}
&\,\,(5,10,29,29,16,27,14)\\[6pt]
\stackrel{\mathrm{Step \;1}}{\longrightarrow} &\left(
\begin{array}{ccccccc}
5&10&29&29&16&27&14\\
1&2&3&4&5&6&7
\end{array}\right)\,\stackrel{\mathrm{Step \;2}}{\longrightarrow}
\left(
\begin{array}{ccccccc}
3&9&27&27&15&27&12\\
1_2&2_1&3_2&4_2&5_1&6_0&7_2
\end{array}\right)\\[6pt]
\stackrel{\mathrm{Step\; 3}}{\longrightarrow}&\left(
\begin{array}{ccccccc}
27&27&27&15&12&9&3\\
3_2&4_2&6_0&5_1&7_2&2_1&1_2
\end{array}\right)\stackrel{\mathrm{Step \;4}}{\longrightarrow}
\left(
\begin{array}{ccccccc}
18&18&18&9&9&6&3\\
3_2&4_2&6_0&5_1&7_2&2_1&1_2
\end{array}\right).
\end{align*}

\section{Labeled Partitions and $q$-Derangements Numbers}\label{Section4}

In this section, we give a combinatorial treatment of a relation on
the  $q$-derangement numbers for $S_n^k$. This relation leads to the
formula of Faliharimalala and Zeng for $d_n^k(q)$.  We will give the
proof for the case $k=3$. It is easy to see that the argument
applies to the general case.

Following Wachs \cite{Wachs89} and Chow \cite{Chow06}, we define the
reduction of a colored permutation $\sigma$ on a set of positive
integers $A=\{a_1<a_2<\cdots <a_k\}$ by substituting the element
$a_i$ with $i$ while keeping the color.  Keep in mind that a positin
$i$ is called a fixed point of a colored permutation $\pi(1)\pi(2)
\cdots \pi(n)$ if $\pi(i)=i_0$. Then the \textit{derangement part}
of a colored permutation $\sigma\in S_n^3$, denoted by $dp(\sigma)$,
is the reduction of the sequence obtained from $\sigma$ by removing
the fixed elements.  For example, $dp(8_0\,1_2\, 5_1\,
4_0\,3_1\,6_0\,7_1\,2_2)=6_0\,1_2\,4_1\,3_1\,5_1\,2_2.$

  Then we have the following extension of the relation due to
  Wachs \cite{Wachs89}:

\begin{theorem}\label{Wachs} Given $\alpha\in
\mathscr{D}_k^3$, for  $0\leq k\leq n$ we have
\begin{align}
\sum_{dp(\sigma)=\alpha,\sigma\in
S_n^3}q^{\mathrm{fmaj}_3(\sigma)}=q^{\mathrm{fmaj}_3(\alpha)}{n\brack
k}_{q^3}.\label{eq generating function of d_n^3}
\end{align}
\end{theorem}

It should be noted that the above theorem can be proved  by the
method of Wachs \cite{Wachs89} which has been extended by Chow
\cite{Chow06} to signed permutations. We will give a combinatorial
proof based on labeled partitions with colored permutations.

For any $\pi=\pi(1)\pi(2)\cdots\pi(k)\in S_k^3$, we can insert a
fixed point $j$ with $1\leq j\leq k+1$ into $\pi$ to obtain a
permutation $\bar{\pi}$ in $S^3_{k+1}$ given by
\begin{align*}
\bar{\pi}=\pi^{\prime}(1)\pi^{\prime}(2)\cdots
\pi^{\prime}(j-1)j_0\pi^{\prime}(j)\cdots\pi^{\prime}(k),
\end{align*}
where
\begin{align*}
\pi^{\prime}(i)=\begin{cases} (c(\pi(i)))d(\pi(i)),\ \ \ \ \
&\text{if}\
d(\pi(i))<j,\\[3pt]
(c(\pi(i)))(d(\pi(i))+1),\ \ \ \ \ &\text{otherwise} .
\end{cases}
\end{align*}

In other words, $\bar{\pi}$ is the unique permutation with $i$ being
a fixed point such that the reduction of the sequence obtained from
$\bar{\pi}$ by deleting the element at position $i$ equals $\pi$.
For example, let $\pi=4_2\,1_0\,2_0\,6_1\,5_1\,3_2$. Then we  get
$5_2\,1_0\,3_0\,2_0\,7_1\,6_1\,4_2$ when we insert $3$ into $\pi$.

\noindent {\it Proof of Theorem \ref{Wachs}}. First, we reformulate
the relation \eqref{eq generating function of d_n^3} in the
equivalent form
\begin{align}
\frac{1}{(q^3;q^3)_n}\sum_{dp(\sigma)=\alpha,\sigma\in
S_n^3}q^{\mathrm{fmaj}_3(\sigma)}
=\frac{1}{(q^3;q^3)_k(q^3;q^3)_{n-k}}q^{\mathrm{fmaj}_3(\alpha)},
\end{align}
and use  labeled partitions to give combinatorial proof of the above
relation. Let $R_{\alpha}$ be the set of colored permutations
$\sigma\in S_n^3$ such that $dp(\sigma)=\alpha$.  We proceed to
establish a bijection $\theta\colon(\lambda,\sigma)\rightarrow
(\beta, \gamma)$ from  $(P_n^3, R_{\alpha})$ to $(P_k^3,P_{n-k}^3)$
such  that
\begin{equation}
|\lambda|+\mathrm{fmaj_3}(\sigma)=|\beta|+|\gamma|+\mathrm{fmaj_3}(\alpha).
\end{equation}
The bijection consists of the following steps.

\noindent
 Step 1. Apply the bijection $g_{\sigma}$ given in
Lemma \ref{lem g_pi} to get a standard labeled partition
$(\lambda^{*}, \sigma)$ from $\lambda$.

\noindent Step 2. Let the fixed points and non-fixed points of
$\sigma$ be $\sigma(i_1),\sigma(i_2),\ldots, \sigma(i_{n-k})$ and
$\sigma(j_1),\sigma(j_2),$ $\ldots,\sigma(j_k)$. We decompose
$\lambda^*$ into two parts, namely,
$\lambda^*(i_1),\lambda^*(i_2),\ldots,\lambda^*(i_{n-k})$ and
$\lambda^*(j_1),\lambda^*(j_2),$ $\ldots,\lambda^*(j_k)$.

Let
$\gamma=(\lambda^*(i_1),\lambda^*(i_2),\ldots,\lambda^*(i_{n-k}))$
and $\beta^*=(\lambda^*(j_1),\lambda^*(j_2),\ldots,
\lambda^*(j_k))$.

\noindent
 Step 3. Apply $g_{\alpha}^{-1}$ to $(\beta^*,\alpha)$
and denote the resulted partition by $\beta$.

To show that the above procedure is feasible, we need to show that
$\beta^*$ generated in Step 2 satisfies the condition that
$(\beta^*,\alpha)$ belongs to $Q_{\alpha}$ so that one can apply
$g_{\alpha}^{-1}$.

Observe that  for any $1\leq q \leq k$,  $\sigma(j_q)$ and
$\alpha(q)$ have the same subscript since $\alpha(q)$ is obtained by
the reduction operation. It follows that
\[\beta^*(q)-\mathrm{c}(\alpha(q))=\lambda^*(j_q)-\mathrm{c}(\alpha(q))\]
is divisible by $3$ for any $1\leq q\leq k$. To prove that
$(\beta^*,\alpha)$ is standard, it suffices to show if
$\sigma(p)>\sigma(q)$ with $\sigma(p+1),\ldots,\sigma(q-1)$ being at
the positions of fixed points, then $\lambda^*_p>\lambda^*_q$. When
$q=p+1$, we conclude that $\lambda^*_p>\lambda^*_q$ from the fact
that $(\lambda^*, \sigma)$ is standard. When $q>p+1$, it is easy to
see that we have either $\sigma(p)>\sigma(p+1)$ or
$\sigma(q-1)>\sigma(q)$. Therefore, we have either
$\lambda^*_p>\lambda^*_{p+1}$ or $\lambda^*_{q-1}>\lambda^*_q$.
Since $\lambda^*$ is a partition, we find that
$\lambda^*_p>\lambda^*_q$. Hence the bijection is well defined.

It remains to show that the above procedure is reversible. We
proceed to construct the inverse map $\eta$ from $(P_k^3,P_{n-k}^3)$
to $(P_n^3, R_{\alpha})$, which consists of three steps.

\noindent
 Step 1. Apply $g_{\alpha}$ to $\beta$ and
denote the resulted partition by $(\tilde{\beta},\alpha)$.

\noindent
 Step 2. Let
$(\tilde{\lambda}^0,\sigma^0)=(\tilde{\beta},\alpha)$. We insert
$\gamma_i$ into $(\tilde{\lambda}^{i-1},\sigma^{i-1})$ to get
$(\tilde{\lambda}^i,\sigma^i)$. Find the first  position $r$ in
$\tilde{\lambda}^{i-1}$ such that the insertion of $\gamma_i$ to
this position will produce a partition. We denote this partition by
$\tilde{\lambda}^i$. Obviously, we have
$\tilde{\lambda}^i_{r-1}>\tilde{\lambda}^i_{r}=\gamma_i$. Suppose
that $\tilde{\lambda}^i_{r}=\cdots
=\tilde{\lambda}^i_{t}>\tilde{\lambda}^i_{t+1}$ for some $t\geq r$.
If $r=t$ then we set $s=r$. Otherwise, from left to right, we look
for a position $s$ satisfying $\sigma^{i-1}(s-1)<s_0\leq
\sigma^{i-1}(s)$ (here we treat $\sigma^{i-1}(r-1)$ as $-\infty$ and
$\sigma^{i-1}(t+1)$ as $\infty$). In this way, we obtain $\sigma^i$
from $\sigma^{i-1}$ by inserting $s_0$ as a fixed point. In fact,
this procedure guarantees that the subsequence
$\sigma^i(r),\sigma^i(r+1),\ldots,\sigma^i(t)$ is increasing. That
is, $(\tilde{\lambda}^i,\sigma^i)$ is a standard labeled partition.
On the other hand, since $\gamma\in P_{n-k}^3$ and each fixed point
has subscript 0, we have $\gamma_i$ is divisible by 3 for each
$1\leq i\leq {n-k}$ and thus $(\tilde{\lambda}^i,\sigma^i)\in
Q_{\sigma^i}$.

\noindent
 Step 3. Apply $g_{\sigma^{n-k}}^{-1}$ to
$(\tilde{\lambda}^{n-k},\sigma^{n-k})$ and  denote the resulted
partition by $\lambda^{n-k}$.

We claim that $\lambda^{n-k}$ and $\sigma^{n-k}$ equal $\lambda$ and
$\sigma$ respectively. Then we see that $\eta$ is the inverse of
$\theta$. From Lemma \ref{lem g_pi}, it is easily seen that
$\beta^*=\tilde{\beta}$. Since $\tilde{\lambda}^{n-k}$ is the
partition obtained from $\tilde{\beta}$ by inserting
$\gamma_1,\ldots, \gamma_{n-k}$, we  have
$\lambda^*=\tilde{\lambda}^{n-k}$.

It is now necessary  to show that $\sigma^{n-k}=\sigma$. It suffices
to verify $\sigma^{n-k}$ and $\sigma$ have the same fixed points. By
removing the  common fixed points, we may assume that the first
fixed point $f$ ($f_0$) of $\sigma$ is different from that of
$f^{\prime}$ ($f^{\prime}_0$) of $\sigma^{n-k}$. We have \[
\sigma(f-1)<f_0\leq \sigma(f+1)-1.\]  Since $f^{\prime}$ is the
first position we aim to find, we have $f^{\prime}<f$. On the other
hand, it is clear that $\lambda^*(f)=\lambda^*(f^{\prime})$. Since
$(\lambda^*,\sigma)$ and $(\lambda^*,\sigma^{n-k})$  are both
standard labeled partitions, we find
\[
\sigma(f^{\prime})<\sigma(f^{\prime}+1)<\cdots<\sigma(f),\] and
\[
\sigma^{n-k}(f^{\prime})<\sigma^{n-k}(f^{\prime}+1)<\cdots<\sigma^{n-k}(f).
\]
Based on the fact that  $\lambda^*(f)=\lambda^*(f^{\prime})$,
$\sigma^{n-k}(f)$ and $\sigma^{n-k}(f^{\prime})$ have the same
subscript,  we conclude that $\sigma^{n-k}(f)$ has the subscript $0$
as $\sigma^{n-k}(f^{\prime})$.

Now we see that $\sigma(f)=f$ and
$\sigma^{n-k}(f^{\prime})=f^{\prime}.$ Since
\[ \sigma(f^{\prime})<\sigma(f^{\prime}+1)<\cdots<\sigma(f)\] and
$\sigma(f)=f$, we can deduce that $\sigma(f^{\prime})\leq
f^{\prime}$. Since $f$ is the first fixed point of $\sigma,$ we
obtain $\alpha(f^{\prime})=\sigma(f^{\prime})<f^{\prime}.$ From the
construction of $\sigma^{n-k},$ if follows that
$\sigma^{n-k}(f^{\prime})\leq \alpha(f^{\prime})<f^{\prime}$ which
contradicts the assumption that
$\sigma^{n-k}(f^{\prime})=f^{\prime}$ ($f^{\prime}$ is a fixed point
of $\sigma^{n-k}$).

Therefore, we have $\sigma=\sigma^{n-k}.$ Again by Lemma \ref{lem
g_pi}, we conclude that  $\lambda=\lambda^{n-k}$. Hence $\eta$ is
the inverse map of $\theta$. This completes the proof. \qed

For example, let $n=8,\lambda=(18,12,12,12,9,9,6,3)$ and
$\sigma=5_2\,1_0\,2_0\,4_0\,8_1\,6_0\,7_1\,3_2$. Then we have
\[g_{\sigma}(\lambda)=\left(
\begin{array}{cccccccc}
29&21&21&21&16&10&5\\
5_2&1_0&2_0&4_0&7_1&6_1&3_2
\end{array}\right).\]
The fixed points of $\sigma$ are $4_0$ and $6_0$, and
$\alpha=dp(\sigma)=4_2\,1_0\,2_0\,6_1\,5_1\,3_2.$ Decomposing
$(29,21,21,21,$ $16,15,10,5)$, we  get $((29,21,21,
16,10,5),(21,15)).$ Applying $g^{-1}_{\alpha}$ to
$\beta^*=(29,21,21, 16,10,5)$ gives $\beta=(18,12,12,9,6,3)$ and
$\gamma=(21,15)$.

Conversely, given $\alpha=4_2\,1_0\,2_0\,6_1\,5_1\,3_2$ and
$(\beta,\gamma)=((18,12,12,9,6,3),(21,15))$, we have
$\tilde{\beta}=(29,21,21, 16,10,5)$. The insertion process is
illustrated as follows:
\begin{align*}\left(
\begin{array}{cccccc}
29&21&21&16&10&5\\
4_2&1_0&2_0&6_1&5_1&3_2
\end{array}\right)&\stackrel{\gamma_1=21}{\longrightarrow}
\left(
\begin{array}{cccccccc}
29&21&21&21&16&10&5\\
5_2&1_0&2_0&4_0&7_1&6_1&3_2
\end{array}\right)\\[6pt] &\stackrel{\gamma_2=15}{\longrightarrow}
\left(
\begin{array}{cccccccc}
29&21&21&21&16&15&10&5\\
5_2&1_0&2_0&4_0&8_1&6_0&7_1&3_2
\end{array}\right).\end{align*}
So we get $\tilde{\lambda}^{n-k}=(29,21,21,21,16,15,10,5)$,
$\sigma^{n-k}=5_2\,1_0\,2_0\,4_0\,8_1\,6_0\,7_1\,3_2$. Finally, we
find $\lambda^{n-k}=g^{-1}_{\sigma^{n-k}}=(18,12,$ $12,12,9,9,6,3).$

\section{Involutions on Labeled Partitions}
\label{Labeled partition and involutions}
\newcommand{\sign}{\mathrm{sign}}

In this section, we give a combinatorial interpretation of the
 formula of Gessel and Simon  in terms of an
involution on labeled partitions. This involution can be easily
extended to  type $B$. Hence we also give a combinatorial proof of a
formula of  Adin, Gessel and Roichman on the signed $q$-counting of
$\mathrm{fmaj}$ indices of signed permutations.

Recall that the sign of a signed permutations is defined in terms of
generators of $B_n$ as a Coxeter group. Consider the
 generating set $\{s_0,s_1,s_2,\ldots,s_{n-1}\}$ of $B_n$,  where
\[s_0:=[-1,2,3,\ldots,n], \quad\textrm{and}\quad
s_i:=[1,2,\ldots,i-1,i+1,i,i+2,\ldots,n]\] for $1\leq i\leq n-1$.
Then the sign of a signed permutation $\pi$ is defined by
\[\sign(\pi):=(-1)^{l(\pi)},\]  where $l(\pi)$ is the
standard length of $\pi$ with respect to the generators of $B_n$.

The following theorem is due to Gessel and Simon \cite{Wachs92}.

\begin{theorem}\label{G-S}
\begin{equation}\sum_{\pi\in
S_n}\mathrm{sign}(\pi)q^{\mathrm{maj}(\pi)}=[1]_{q}[2]_{-q}[3]_{q}[4]_{-q}\cdots[n]_{(-1)^{n-1}q}.\label{G-S
formula}\end{equation}
\end{theorem}

A combinatorial proof of the above formula has been given by Wachs
\cite{Wachs92}. Here we will give an involution on labeled
partitions and we will show that this involution can be easily
extended to the following type $B$ formula due to  Adin, Gessel and
Roichman \cite{Adin05}.

\begin{theorem}\label{A-G-R}
\begin{equation}\sum_{\pi\in
B_n}\mathrm{sign}(\pi)q^{\mathrm{fmaj}(\pi)}=[2]_{-q}[4]_{q}\cdots[2n]_{(-1)^{n}q}.\label{A-G-R-F}
\end{equation}
\end{theorem}

To describe our involution on labeled partitions as a proof of the
formula \eqref{G-S formula}, we may reformulate in
 an equivalent form:
\begin{equation}
\frac{\sum_{\pi\in
S_n}\mathrm{sign}(\pi)q^{\mathrm{maj}(\pi)}}{(q;q)_n}=\frac{1}{(1-q)(1+q)(1-q)(1+q)\cdots(1-(-1)^{n-1}q)}.
\label{G-S-E}
\end{equation}

 \noindent{\it Proof of Theorem \ref{G-S}.} We consider two cases
 according to the parity of $n$.

Case 1. $n$ is even, i.e., $n=2k.$ Then (\ref{G-S-E}) takes the form
\begin{equation}
\frac{\sum_{\pi\in
S_{2k}}\mathrm{sign}(\pi)q^{\mathrm{maj}(\pi)}}{(q;q)_{2k}}=\frac{1}{(1-q^2)^k}.
\label{G-S-e}
\end{equation}
Clearly, the right hand side of (\ref{G-S-e}) is the generating
function of sequences $(a_1,a_2,\ldots,a_{2k-1},a_{2k})$ satisfying
$a_{2i-1}=a_{2i}$ for $i=1,2,\ldots,k$. It is also easy to see that
the left hand side of (\ref{G-S-e}) is the generating function of
labeled partitions  on $S_n$ with at most $2k$ parts under the
assumption that a labeled partition $(\lambda,\pi)$ carries the sign
of the permutation $\pi$. To be more specific,  such labeled
partitions are called signed labeled partitions. We proceed to
construct an involution on the set $H$ of  signed labeled partitions
$(\lambda, \pi)$ such that the generating function of the fixed
points of this involution equals the right hand side of
(\ref{G-S-e}). This involution consists of three steps.

\noindent
 Step 1. Let $(\lambda, \pi)$ be a labeled partition such that $\pi\in S_{2k}$ and
$\lambda=(\lambda_1,\lambda_2,\ldots,\lambda_{2k})$ with
$\lambda_1\geq \lambda_2\geq \cdots\geq \lambda_{2k}\geq 0$. If
$|\pi^{-1}(1)-\pi^{-1}(2)|\neq 1$, then we define
$$\phi^1(\pi)(i)=\left\{
\begin{array}{cl}
\pi(i),&i\neq \pi^{-1}(1) \textrm{ and }\pi^{-1}(2),\\[3pt]
2,&i=\pi^{-1}(1),\\[3pt]
1,&i=\pi^{-1}(2).\\
\end{array}\right.$$
Obviously,  $(\lambda,\pi)$ and $(\lambda,\phi^1(\pi))$  have
opposite signs and $\mathrm{maj}(\pi)=\mathrm{maj}(\phi^1(\pi))$.
Therefore, we have
$\mathrm{maj}(\pi)+|\lambda|=\mathrm{maj}(\phi^1(\pi))+|\lambda|$,
and so these two elements cancel each other. If
$|\pi^{-1}(1)-\pi^{-1}(2)|=1$, then we see that
$\mathrm{maj}(\pi)\neq\mathrm{maj}(\phi^1(\pi))$.

We now use $H^1$ to denote the set of  signed labeled partitions
$(\lambda,\pi)$ such that $|\pi^{-1}(1)-\pi^{-1}(2)|= 1$. Repeating
the above procedure, we continue to cancel out some elements in
$H^1$. At this time,  we consider the positions of the elements $3$
and $4$. Similarly, if $|\pi^{-1}(3)-\pi^{-1}(4)|\neq 1$, then we
define
$$\phi^2(\pi)(i)=\left\{
\begin{array}{cl}
\pi(i),&i\neq \pi^{-1}(3)\textrm{ and }\pi^{-1}(4),\\[3pt]
4,&i=\pi^{-1}(3),\\[3pt]
3,&i=\pi^{-1}(4).\\
\end{array}\right.$$
It follows that $(\lambda,\pi)$ and $(\lambda,\phi^2(\pi))$  have
the opposite signs and that
\[ \mathrm{maj}(\pi)+|\lambda|=\mathrm{maj}(\phi^2(\pi))+|\lambda|.\]
In other words, the two elements cancel out  in the set $H^1$.

Now, we use $H^2$ to denote the subset of $H^1$ such that
$|\pi^{-1}(3)-\pi^{-1}(4)|=1$. Iterating this process, we may
consider the elements $\{5,6\},\{7,8\},\ldots,\{2k-1,2k\}$ and
denote the set obtained at the last step by $H^k$. Finally, we
obtain $H^k\subseteq H^{k-1}\subseteq\cdots \subseteq H^1$. In the
intermediate steps, we can defined the functions $\phi^i$ for $i=1,
2, \ldots, k$. It is not difficult to see that the labeled partition
$(\lambda, \pi)$ in $H^k$ has the property that
\[ |\pi^{-1}(1)-\pi^{-1}(2)|= 1, |\pi^{-1}(3)-\pi^{-1}(4)|=
1,\ldots, |\pi^{-1}(2k-1)-\pi^{-1}(2k)|= 1.\] Namely, any odd number
$2i-1$ is next to $2i$ in $\pi$ for all $i=1,\ldots, k$.

\noindent Step 2. For any labeled partition
\begin{align*}(\lambda,\pi)=\left(
\begin{array}{cccccc}
\lambda_1&\cdots&\lambda_{\pi^{-1}(2)}&\lambda_{\pi^{-1}(1)}&\cdots&\lambda_{2k}\\
\pi(1)&\cdots&2&1&\cdots&\pi(2k)
\end{array}\right),
\end{align*}
we define $(f^1(\lambda),g^1(\pi))$  to be the labeled partition
\begin{align*}
(f^1(\lambda),g^1(\pi))=\left(
\begin{array}{cccccc}
\lambda_1+1&\cdots&\lambda_{\pi^{-1}(2)}+1&\lambda_{\pi^{-1}(1)}&\cdots&\lambda_{2k}\\
\pi(1)&\cdots&1&2&\cdots&\pi(2k)
\end{array}\right),
\end{align*}
where $f^1(\lambda)$ is the partition obtained from $\lambda$ by
adding $1$ to the first $\pi^{-1}(2)$ parts of $\lambda$ and
$g^1(\pi)$ is the permutation obtained from $\pi$ by exchanging the
positions of $1$ and $2$.

 Clearly, $(\lambda,\pi)$ and $(f^1(\lambda),g^1(\pi))$ have
opposite signs. Also, we have
\[\mathrm{maj}(\pi)+|\lambda|=\mathrm{maj}(g^1(\pi))+|f^1(\lambda)|.\]
 Therefore  $(\lambda,\pi)$ and
$(f^1(\lambda),g^1(\pi))$ cancel out in $H^k$. Notice that the
resulted labeled partition $(f^1(\lambda),g^1(\pi))$ has the
additional property that $f^1(\lambda)_{\pi^{-1}(1)}$ is greater
than $f^1(\lambda)_{\pi^{-1}(2)}$.  By inspection, we
 see that after cancellation, the remaining elements in  $H^k$ are of the following form
 \begin{align*}(\lambda,\pi)=\left(
\begin{array}{cccccc}
\lambda_1&\cdots&\lambda_{\pi^{-1}(1)}&\lambda_{\pi^{-1}(2)}&\cdots&\lambda_{2k}\\
\pi(1)&\cdots&1&2&\cdots&\pi(2k)
\end{array}\right)
\end{align*}
where $\lambda_{\pi^{-1}(1)}=\lambda_{\pi^{-1}(2)}$. Let $H^k_1$
denote the set of remaining elements in $H^k$ that of the above
form.

We continue the above process for the $H_1^k$ with respect the
relative positions of $3$ and $4$. It is easy to check that for any
labeled partition $(\lambda,\pi)$ in $H^k_1$, 1 appears before 2 in
$\pi$ and $\lambda_{\pi^{-1}(1)}=\lambda_{\pi^{-1}(2)}$. Now, for
any element $(\lambda,\pi)\in H^k_1$, if
\begin{align*}(\lambda,\pi)=\left(
\begin{array}{cccccc}
\lambda_1&\cdots&\lambda_{\pi^{-1}(4)}&\lambda_{\pi^{-1}(3)}&\cdots&\lambda_{2k}\\
\pi(1)&\cdots&4&3&\cdots&\pi(2k)
\end{array}\right),
\end{align*}
then we can find  another labeled partition $(f^2(\lambda),g^2(\pi))
\in H^k_1$
\begin{align*}
(f^2(\lambda),g^2(\pi))=\left(
\begin{array}{cccccc}
\lambda_1+1&\cdots&\lambda_{\pi^{-1}(4)}+1&\lambda_{\pi^{-1}(3)}&\cdots&\lambda_{2k}\\
\pi(1)&\cdots&3&4&\cdots&\pi(2k)
\end{array}\right).
\end{align*}
Again,  $(\lambda,\pi)$ and $(f^2(\lambda),g^2(\pi))$ cancel each
other in $H^k_1$. Notice that  $f^2(\lambda)_{\pi^{-1}(3)}$ is
greater than $f^2(\lambda)_{\pi^{-1}(4)}$.  So the remaining labeled
partitions after the above cancelation are of the form
\begin{align*}(\lambda,\pi)=\left(
\begin{array}{cccccc}
\lambda_1&\cdots&\lambda_{\pi^{-1}(3)}&\lambda_{\pi^{-1}(4)}&\cdots&\lambda_{2k}\\
\pi(1)&\cdots&3&4&\cdots&\pi(2k)
\end{array}\right),
\end{align*}
where $\lambda_{\pi^{-1}(3)}=\lambda_{\pi^{-1}(4)}$. Then we can
denote the set of the remaining labeled partitions by $H^k_2$ and
continue the above process.
 Eventually, we get  $H^k_k\subseteq
H^{k}_{k-1}\subseteq\cdots\subseteq H^k_1$. Moreover, in the process
we have defined the functions $f^i$ and $g^i$ for $i=1, 2, \ldots,
k$.

It is easy to see that for any labeled partition  $(\lambda, \pi)$
in $H^k_k$ and for any $i\in \{1,\ldots, k\}$, $2i-1$ appears
immediately before $2i$  and
$\lambda_{\pi^{-1}(2i-1)}=\lambda_{\pi^{-1}(2i)}$. Clearly, all the
labeled partitions in $H_k^k$ have positive signs.

\noindent Step 3. Permute the columns of the labeled partitions
$(\lambda, \pi)$ in $H^k_k$ so that the elements in $\pi$ are
rearranged in increasing order. Taking the first row of the resulted
two row array, we will get a sequence
$(a_1,a_2,\ldots,a_{2k-1},a_{2k})$ such that
$a_{2i-1}=a_{2i}\,(i=1,\ldots, k)$ whose generating function is the
right hand side of \eqref{G-S-e}.

It is easy to  see that  the relation \eqref{G-S-e} can be justified
by the above algorithm. Hence Theorem \ref{G-S} holds when $n$ is
even.

\medskip

\noindent Case 2. $n$ is odd, i.e., $n=2k+1.$ We need to show that
\begin{equation}
\frac{\sum_{\pi\in
S_{2k+1}}\mathrm{sign}(\pi)q^{\mathrm{maj}(\pi)}}{(q;q)_{2k+1}}=\frac{1}{(1-q^2)^k(1-q)}.\label{G-S
for n is odd}
\end{equation}

This case is analogous to the case when $n$ is even. We may employ
the same operations in Step 1 and Step 2 by ignoring the element
$2k+1$ while making the pairs $\{1, 2,\}, \{3, 4\}, \ldots, \{2k-1,
2k\}$. The only difference lies in Step 3 when we take the first row
of the resulted two row array, we  get a sequence
$(a_1,a_2,\ldots,a_{2k-1},a_{2k},a_{2k+1})$ such that
$a_{2i-1}=a_{2i}\,(i=1,\ldots, k)$. Moreover, $a_{2k+1}$ can be any
positive integer. This completes the proof of the relation
\eqref{G-S for n is odd}.

In fact, we have constructed a sign reversing involution
\[ (\theta,\chi)\colon (\lambda, \pi) \rightarrow
(\theta(\lambda),\chi(\pi)).\]Specifically, the map $(\theta, \chi)$
is defined by
$$(\theta(\lambda),\chi(\pi))=\left\{
\begin{array}{lll}
&(\lambda,\phi^1(\pi)),&\text{if}(\lambda,\pi)\in H\setminus
 H^1,\\[5pt]
&(\lambda,\phi^2(\pi)),&\text{if}(\lambda,\pi)\in H^1\setminus
 H^2,\\[5pt]
&\cdots&\\[5pt]
&(\lambda,\phi^k(\pi)),&\text{if}(\lambda,\pi)\in
H^{k-1}\setminus H^k,\\[5pt]
&(f^1(\lambda),g^1(\pi)),&\text{if}(\lambda,\pi)\in
H^k\setminus H^k_1,\\[5pt]
&(f^2(\lambda),g^2(\pi)),&\text{if}(\lambda,\pi)\in
H^k_1\setminus H^k_2,\\[5pt]
&\cdots&\\[5pt]
&(f^k(\lambda),g^k(\pi)),&\text{if}(\lambda,\pi)\in
H^k_{k-1}\setminus
 H^k_k,\\[5pt]
&(\lambda,\pi),&\text{if}(\lambda,\pi)\in H^k_k,
\end{array}\right.$$
where $\phi^i(\pi)$, $f^i(\lambda)$ and $g^i(\pi)$ are defined in
the above algorithm. It is easy to verify that the map induces sign
reversing, that is, if $(\lambda, \pi)$ is not a fixed point of the
map $(\theta, \chi)$, then we have
$\mathrm{sign}(\theta(\lambda),\chi(\pi))=-\mathrm{sign}(\lambda,\pi)$
and
$|\theta(\lambda)|+\mathrm{maj}(\chi(\pi))=|\lambda|+\mathrm{maj}(\pi)$.
The fixed points of the map $(\theta, \chi)$ correspond to the right
hand side of \eqref{G-S formula}. This completes the proof. \qed

We now turn to the Theorem \ref{A-G-R}, and we need a
characterization of the length function of signed permutations
\cite[Propostion 3.1 and Corollary 3.2]{Bre}.

\begin{lem}\label{lengthlemma}
Let $\sigma\in B_n$, we have
\[l(\sigma)=\textrm{\inv}\,(\sigma)+\sum_{\{1\leq i\leq n|\sigma(i)<0\}}|\sigma(i)|,\]
where $\textrm{inv}\,(\sigma)$ is defined  with respect to the order
\[\bar{n}<\cdots<\bar{1}<1<\cdots<n.\]
\end{lem}

Note that in the definition of the $\mathrm{fmaj}$ on $B_n$ we have
imposed the order
\[\bar{1}<\cdots<\bar{n}<1<\cdots<n \]
or in the notation of colored permutations,
\[1_1<\cdots<n_1<1_0<\cdots<n_0.\]

The above lemma is useful for the construction of a sign reversing
involution for the formula (\ref{A-G-R-F}) for $B_n$. Given a signed
permutation $\sigma\in B_n$, we may construct a signed permutation
$\sigma'$ as follows. If $1$ and $2$ have different signs or $1$ and
$2$ have the same sign but are not adjacent in $\sigma$, then we
exchange $1$ and $2$ without changing the signs.  By Lemma
\ref{lengthlemma}, we see that the  $\sigma'$ and $\sigma$ have
opposite signs and
$\mathrm{fmaj}(\sigma)=\mathrm{fmaj}(\sigma^{\prime})$.

For example, let $\sigma=4_0\,2_1\,5_1\,1_0\,3_1$. Then we have
$\sigma'=4_0\,1_1\,5_1\,2_0\,3_1$. Clearly,  $\sigma$ and $\sigma'$
have opposite signs.

Using the above sign change rule, we can  extend the involution for
Theorem \ref{G-S} to Theorem \ref{A-G-R}.  The detailed proof is
omitted. \qed

\vspace{0.5cm}
 \noindent{\bf Acknowledgments.}  This work was supported by  the 973
Project, the PCSIRT Project of the Ministry of Education, the
Ministry of Science and Technology, and the National Science
Foundation of China.

\end{document}